\documentclass[11pt]{amsart}
\usepackage{amssymb,amsmath,amsthm}
\usepackage{graphicx}

\usepackage{a4wide}

\numberwithin{equation}{subsection}

\newtheorem{thm}{Theorem}[subsection]
\newtheorem{prop}[thm]{Proposition}
\newtheorem{lem}[thm]{Lemma}
\newtheorem{cor}[thm]{Corollary}
\newenvironment{rem}{\vspace{3mm}\noindent
{\bf Remark.}}{\vspace{3mm}}

\newenvironment{rems}{\vspace{3mm}
\noindent {\bf Remarks.}}{\vspace{3mm}}

\newcommand{\Pf}{\noindent {\it Proof}}

\newcommand*{\Dd}{\mathop{\mathrm D\kern0pt}\nolimits}

\long\def\comment#1{}

\newcommand{\Tr}{\operatorname{Tr}}

\newcommand{\OO}{{\mathcal O}}

\newcommand{\SL}{\operatorname{SL}}

\newcommand{\lan}{\langle}
\newcommand{\ran}{\rangle}

\newcommand{\de}{\delta}
\newcommand{\eps}{\epsilon}
\renewcommand{\ker}{\operatorname{ker}}

\newcommand{\id}{\operatorname{id}}

\newcommand{\ov}{\overline}
\newcommand{\we}{\wedge}
\renewcommand{\Im}{\operatorname{Im}}

\newcommand{\EE}{{\mathcal E}}

\newcommand{\Om}{\Omega}
\newcommand{\dbar}{\overline{\partial}}
\newcommand{\pa}{\partial}
\newcommand{\Hom}{\operatorname{Hom}}

\newcommand{\Ext}{\operatorname{Ext}}

\renewcommand{\a}{\alpha}
\renewcommand{\b}{\beta}
\newcommand{\om}{\omega}

\newcommand{\la}{\lambda}
\renewcommand{\th}{\theta}
\newcommand{\C}{{\Bbb C}}

\newcommand{\R}{{\Bbb R}}
\newcommand{\Z}{{\Bbb Z}}
\newcommand{\Q}{{\Bbb Q}}

\newcommand{\ot}{\otimes}

\newcommand{\sub}{\subset}
\newcommand{\ed}{\qed\vspace{3mm}}

\newcommand{\join}{\operatorname{join}}

\begin{document}

\title{$A_{\infty}$-algebra of an elliptic curve and Eisenstein series}
\author{Alexander Polishchuk}
\address{Department of Mathematics, University of Oregon, Eugene, OR 97403}  
\email{apolish@uoregon.edu}
\thanks{Supported in part by NSF grant}


\begin{abstract} We compute explicitly the $A_{\infty}$-structure on the algebra $\Ext^*(\OO_C\oplus L,\OO_C\oplus L)$,
where $L$ is a line bundle of degree $1$ on an elliptic curve $C$. The answer involves higher derivatives of
Eisenstein series.
\end{abstract}
\maketitle

\section*{Introduction}

The bounded derived category $D^b(X)$ of coherent sheaves is an important invariant of an algebraic variety $X$
(see \cite{BO} for a survey). This category can be described in a purely algebraic
way if one has a {\it generator}, i.e., an object $G$ generating $D^b(X)$ as a triangulated category.
With such an object one can associate a graded algebra $E_G=\oplus_{n\in\Z}\Ext^n(G,G)$. However,
in order to recover $D^b(X)$ from $E_G$ one has to take into account certain higher products which fit together
into a structure of an {\it $A_{\infty}$-algebra} on $E_G$ (see \cite{Keller-intro} for an introduction into
$A_{\infty}$-algebras). Namely, 
one can realize $E_G$ as the cohomology of a dg-algebra and then apply a general algebraic construction
that gives an $A_{\infty}$-structure on such cohomology (see \cite{Merk}). 
This $A_{\infty}$-structure is minimal in the sense
that $m_1=0$, and is canonical up to $A_{\infty}$-equivalence.
Now the category $D^b(X)$ can be shown to be
equivalent to the derived category of perfect $A_{\infty}$-modules over $A$ (see \cite{Keller}, Thm. 3.1 or
\cite{Lef}, 7.6). 

Thus, it is of interest to compute explicitly higher products on algebras of the form $\Ext^*(G,G)$ as above.
In this paper we solve this problem in the case when $X$ is a complex elliptic curve and $G=\OO_X\oplus L$,
where $L$ is a line bundle of degree $1$. Namely, we compute the $A_{\infty}$-structure arising 
from the harmonic representatives (with respect to natural metrics) in the Dolbeault complex computing 
$\Ext^*(G,G)$. The resulting formulas involve higher derivatives of Eisenstein series
(see Theorem \ref{main-thm}). More precisely, we have to use the well-known 
non-holomorphic (but modular) modification $e^*_2$ of the standard Eisenstein series $e_2$ along
with all the higher Eisenstein series $e_{2k}$, $k\ge 2$. 

It is interesting that Eisenstein series appear not in their usual form but rather as some rapidly decreasing 
series, similar to those considered in \cite{P-Wz} (see Theorem \ref{eis-thm}).
The $A_{\infty}$-constraint gives rise to some quadratic relations involving derivatives of Eisenstein series,
some of them well known (see Proposition \ref{id-prop}).

\noindent{\it Acknowledgment.} I am grateful to Alexandr Usnich and Dmytro Shklyarov for helpful
discussions. This paper was written during the stay at the
Institut des Hautes \'Etudes Scientifiques. I'd like to thank this institution for hospitality and support.

\section{Eisenstein series}

\subsection{Definitions}
Let us recall some basic definitions and facts (see \cite{Weil}, ch.~III).
We consider $C^{\infty}$-functions $F(\om_1,\om_2)$ on the space of all oriented bases of $\C=\R^2$.
Recall that such a function is said to be of weight $k\in\Z$ if
$$F(\la\om_1,\la\om_2)=\la^{-k}F(\om_1,\om_2).$$
Such $F$ is called modular (with respect to $\SL(2,\Z)$) if it is invariant with respect
to $\SL(2,\Z)$ base changes of $(\om_1,\om_2)$ and $F(1,\tau)=f(e^{2\pi i\tau})$, where $f(q)$ is meromorphic
at $q=0$. The Eisenstein series $e_{2k}$ for $k\ge 2$ is defined by
$$e_{2k}(\om_1,\om_2)=\sum_{\om\in L\setminus\{0\}}\frac{1}{\om^{2k}}$$
where $L=\Z\om_1+\Z\om_2$. The function $e_{2k}$ is modular of weight $2k$.
One can also consider the analogous series for $k=1$ using Eisenstein's summation rule:
$$e_2(\om_1,\om_2)=\sum_m\sum_{n;n\neq 0\text{ if }m=0}\frac{1}{(m\om_2+n\om_1)^2}.$$
The function $e_2$ is not modular but admits a simple non-holomorphic correction that makes it $\SL(2,\Z)$-invariant.
Namely, let us set
$$e^*_2(\om_1,\om_2)=e_2(\om_1,\om_2)-\frac{\pi}{a(L)}\frac{\bar{\om_1}}{\om_1},$$
where $a(L)=\Im(\ov{\om}_1\om_2)$ is the area of $\C/L$.
Then $e^*_2(\om_1,\om_2)$ is $\SL(2,\Z)$-invariant (and of weight $2$). 
For convenience we set also $e^*_{2k}=e_{2k}$ for $k\ge 2$.
Eisenstein series appear as coefficients in the expansion of the
Weierstrass zeta-function (cf. \cite{Weil}, ch.~III, formula (9)):
\begin{equation}\label{zeta-exp}
\zeta(z;\om_1,\om_2)=\frac{1}{z}-\sum_{k\ge 2}e_{2k}(\om_1,\om_2)z^{2k-1}.
\end{equation}

Following \cite{Katz}, sec. 1.5, we consider the normalized Weil operator 
$$W=-\frac{\pi}{a(L)}(\ov{\om}_1\frac{\pa}{\pa\om_1}+\ov{\om_2}\frac{\pa}{\pa\om_2}).$$
This operator is $\SL(2,\Z)$-invariant and is of weight two.
Slightly modifying the definition in \cite{Weil}, ch.VI, for a pair of integers $b>a\ge 0$ of different parity we set
$$g_{a,b}=g_{b,a}=(b-a)!W^a(e^*_{b-a+1}).$$
Note that $(a(L)/\pi)^a\cdot g_{a,b}$ differs by a rational factor from Weil's $e^*_{a,b+1}$.
As shown in \cite{Weil}, sec.~VI.5, $g_{a,b}$ is a polynomial in $e^*_2,e_4,\ldots,e_{a+b+1}$ with
rational coefficients.

\subsection{Presentation by rapidly decreasing series}
For $m,n\in\Z$ and a lattice $L\sub\C$ we set
$$f_{m,n}(L)=(\frac{\pi}{a(L)})^m\cdot \sum_{\om\in L\setminus\{0\}}\frac{\bar{\om}^m}{\om^n}\exp(-\frac{\pi}{a(L)}|\om|^2).$$
Note that $f_{m,n}=0$ unless $m+n$ is even (as one can see making the substitution $\om\mapsto-\om$).

\begin{thm}\label{eis-thm}
For $n=2k\ge 2$ one has
$$e^*_{n}=\frac{2}{(n-1)!}f_{n-1,1}+\sum_{m=2}^{n}\frac{1}{(n-m)!}f_{n-m,m}=
\sum_{\om\neq 0}\frac{P_{n-1}(\frac{\pi}{a(L)}|\om|^2)}{\om^n}\exp(-\frac{\pi}{a(L)}|\om|^2),$$
where 
$$P_{n-1}(z)=\frac{2z^{n-1}}{(n-1)!}+\sum_{m=2}^n\frac{z^{n-m}}{(n-m)!}.$$
\end{thm}

\Pf . This follows easily from Theorem 1 of \cite{P-Wz} stating that
$$\zeta(z;\om_1,\om_2)-z_1\eta_1-z_2\eta_2=
\sum_{\om\in L}\frac{\exp(-\frac{\pi}{a(L)}|\om+z|^2)}{\om+z}-
\sum_{\om\in L\setminus\{0\}}\frac{\exp(-\frac{\pi}{a(L)}|\om|^2+2\pi i\frac{\Im(\ov{\om}z)}{a(L)})}{\om},$$
where $z=z_1\om_1+z_2\om_2$ with $z_1,z_2\in\R$.
Indeed, we can use the expansion \eqref{zeta-exp} to check the assertion for $n\ge 4$: one has to subtract
$1/z$ from both parts of the above identity, apply $(\frac{\pa}{\pa z})^{n-1}$ and then evaluate at $z=0$.
The case $n=2$ is slightly different: we again subtract $1/z$ from both parts, then apply
$\frac{\pa}{\pa z_1}$ (taking $(z_1,z_2)$ as independent variables) and evaluate at $z=0$.
The required formula follows the fact that $e_2=-\eta_1/\om_1$ (see e.g., \cite{Katz}, sec. 1.2).
\ed

\begin{rem} The fact that $e^*_{n}$ is holomorphic in $(\om_1,\om_2)$ for $n>2$ is equivalent to the identity 
$$2f_{n-1,-1}=(n-1)f_{n-2,0}.$$
For $n=2$ we have instead $2f_{1,-1}=f_{0,0}+1$. These identities can be derived from the Poisson summation
formula and Fourier self-duality of $\exp(-\frac{\pi}{a(L)}|z|^2)$ (see \cite{P-Wz}, sec. 1.1, Remark 1; for $n>2$ one also has to use differentiation). 
\end{rem}

One can immediately check that 
$$W(f_{m,n})=f_{m+2,n}+nf_{m+1,n+1}.$$
Hence, from Theorem \ref{eis-thm} we get the following formula for $g_{a,b}$.

\begin{cor}\label{main-cor} 
For a pair of integers $a,b\ge 0$ of different parity one has
$$g_{a,b}=\sum_{k\ge 0}k!\left({a\choose k}+{b\choose k}\right)f_{a+b-k,k+1}$$
\end{cor}

\section{Minimal $A_{\infty}$-algebra of an elliptic curve}

\subsection{General construction}

Let us first recall the general construction of the $A_{\infty}$-structure
on the cohomology of a dg-algebra $(A,d)$ equipped with a projector
$\Pi:A\to B$ onto a subspace of $\ker(d)$ and a homotopy operator $Q$
such that $1-\Pi=dQ+Qd$. 
Merkulov's formula for this $A_{\infty}$-structure (see \cite{Merk}) was rewritten in
\cite{KS} as a sum over trees:
$$m_n(b_1,\ldots,b_n)=-\sum_{T}\eps(T)m_T(b_1,\ldots,b_n).$$
Here $T$ runs over all oriented planar rooted $3$-valent trees with $n$ leaves (different from the root)
marked by $b_1,\ldots,b_n$
left to right, and the root marked by $\Pi$ (we draw the tree in such a way that leaves are above, and
every vertex has two edges coming from above and one from below). The expression
$m_T(b_1,\ldots,b_n)$ is obtained by going down from leaves to the root, applying multiplication in $A$ at
every vertex and applying the operator $Q$ at every inner edge (see \cite{KS}, sec. 6.4, for details). 
The sign $\eps(T)$ has form
$$\eps(T)=\prod_v (-1)^{|e_1(v)|+(|e_2(v)|-1)\deg(e_1(v))},$$
where $v$ runs through vertices of $T$, $(e_1(v),e_2(v))$ is the pair of edges above $v$,
for an edge $e$ we denote by $|e|$ the total number of leaves above $e$ and by $\deg(e)$ the sum of
degrees of all leaves above $e$ (recall that leaves are marked by $b_i$).

\begin{lem}\label{unit-lem} 
Assume in addition that $\Pi Q=Q\Pi=Q^2=0$. Let $(b_1,\ldots,b_n)$ be a collection of elements in $B$, where $n\ge 3$,
such that $b_i=1$ for some $i$. Then $m_n(b_1,\ldots,b_n)=0$.
\end{lem}

\Pf . It is convenient to use Merkulov's original formula 
$$m_n(b_1,\ldots,b_n)=\Pi\la_n(b_1,\ldots,b_n),$$
where $\la_n:A^{\ot n}\to A$ are defined for $n\ge 2$ by the following recursion: $\la_2(a_1,a_2)=a_1a_2$,
\begin{align*}
&\la_n(a_1,\ldots,a_n)=\pm Q(\la_{n-1}(a_1,\ldots,a_{n-1}))\cdot a_n\pm a_1\cdot Q(\la_{n-1}(a_2,\ldots,a_n))+\\
&\sum_{k+l=n;k,l\ge 2}\pm Q(\la_k(a_1,\ldots,a_k))\cdot Q(\la_l(a_{k+1},\ldots,a_n)).
\end{align*}
Since, $\Pi Q=0$, it is enough to prove that $\la_n(b_1,\ldots,b_n)\in Q(A)$. Let us use induction in $n$.
In the case $n=3$ we have
$$\la_3(b_1,b_2,b_3)=Q(b_1b_2)b_3\pm b_1Q(b_2b_3)$$
and the assertion follows immediately from the fact that $Q(B)=0$.  
Suppose now that $n\ge 4$ and the assertion holds for all $n'<n$. Since $Q^2=0$, 
the induction assumption easily implies that the first two terms in the recursive formula for $\la_n$ belong to $Q(A)$.
Similarly, all the remaining terms vanish if $n\ge 5$. In the case $n=4$ the term $Q(b_1b_2)\cdot Q(b_3b_4)$ also vanishes because either $b_1b_2\in B$ or $b_3b_4\in B$ and $Q(B)=0$.
\ed

\subsection{The case of an elliptic curve}

Let $C=\C/(\Z\oplus\Z\tau)$ be a complex elliptic curve. We denote by
$L$ the holomorphic line bundle of degree $1$ on $C$, such that
the theta-function $\th(z,\tau)$ descends to a global section of $L$. We consider the Dolbeault dg-algebra
$$A=(\Om^{0,*}\ot {\EE}nd(\OO_C\oplus L),\dbar).$$
Its cohomology $B$ is the direct sum of the following components:

\noindent
(i) $\Hom(\OO,\OO)$ and $\Hom(L,L)$, both generated by identity maps;

\noindent
(ii) $\Hom(\OO,L)$, one-dimensional space;

\noindent
(iii) $\Ext^1(L,\OO)$, one-dimensional space;

\noindent
(iv) $\Ext^1(\OO,\OO)$ and $\Ext^1(L,L)$, both isomorphic to the one-dimensional space $H^1(\OO)$.

To construct the homotopy operator $Q$, as in \cite{P-hi-pr}, we use 
the flat metric on $C$ and the hermitian metric on $L$ given by
$$(f,g)=\int_C f(z)\ov{g(z)}\exp(-2\pi \frac{\Im(z)^2}{\Im(\tau)}) dxdy,$$
where $z=x+iy$. Then we set $Q=\dbar^*G$, where $G$ is the Green operator corresponding to the Laplacian
$\dbar^*\dbar+\dbar\dbar^*$. Then $B\sub A$ is exactly the space of harmonic forms, and $\Pi:A\to B$
is the orthogonal projection.
Let us fix the following harmonic generators in the above components:

\noindent
(i) $\id_{\OO}$ and $\id_L$;

\noindent
(ii) $\th=\th(z,\tau)$ viewed as a holomorphic section of $L$;

\noindent
(iii) $\eta:=\sqrt{2\Im(\tau)}\cdot\ov{\th(z,\tau)}\exp(-2\pi\frac{\Im(z)}{\Im(\tau)^2})d\ov{z}$ viewed as a $(0,1)$-form with values in $L^{-1}$;

\noindent
(iv) $\xi=d\ov{z}$. When it is viewed as an element of $\Ext^1(L,L)$ we write $\xi_L$.

Note that we have a natural symmetric bilinear pairing on $A=A^0\oplus A^1$ given by 
$$\lan \a,\b\ran=\frac{1}{2i\Im(\tau)}\cdot\int_C \Tr(\a\circ \b)\we dz,$$
where $\a$ and $\b$ are homogeneous elements such that $\deg(\a)+\deg(\b)=1$.
The normalization is chosen in such a way that
$$\lan \xi,1\ran=\frac{1}{2i\Im(\tau)}\cdot\int_C d\ov{z}\we dz=1.$$
By Serre duality, the induced pairing between $B^0$ and $B^1$ is nondegenerate.
Also, by Theorem 1.1 of \cite{P-hi-pr}, the $A_{\infty}$-structure on $B$ satisfies
the following cyclic symmetry: 
\begin{equation}\label{cyclic-sym}
\lan m_n(\a_1,\ldots, \a_n),\a_{n+1}\ran=(-1)^{n(\deg(\a_1)+1)}\lan\a_1,m_n(\a_2,\ldots,\a_{n+1})\ran.
\end{equation}

The product $m_2$ on $B$ is just the induced product on cohomology. The only interesting products
are $m_2(\th,\eta)\in\Ext^1(\OO,\OO)$ and $m_2(\eta,\th)\in\Ext^1(L,L)$. Both are proportional to the generator $\xi$.
To find the coefficient of proportionality it is enough to compute
$$\lan m_2(\th,\eta),\id_{\OO}\ran=\lan m_2(\eta,\th),\id_L\ran=\frac{1}{2i\Im(\tau)}\cdot\int_C \th\cdot\eta\we dz.$$
The above integral is well known:
$$\frac{\sqrt{2\Im(\tau)}}{2i\Im(\tau)}\cdot\int_C \th(z,\tau)\ov{\th(z,\tau)}\exp(-2\pi\frac{\Im(z)}{\Im(\tau)^2})d\ov{z}\we dz=1.$$
Thus, we obtain
\begin{equation}\label{m2-for}
m_2(\th,\eta)=\xi, \ \ m_2(\eta,\th)=\xi_L.
\end{equation}

By Lemma \ref{unit-lem}, every higher product $m_n$ containing $\id_{\OO}$ or $\id_L$ vanishes.
Together with the cyclic symmetry \eqref{cyclic-sym} this implies that the only potentially nonzero higher products are of 
the following types:

\noindent
(I) $m_n((\xi)^{a},\th,(\xi_L)^b,\eta,(\xi)^c,\th,(\xi_L)^d)\in\Hom(\OO,L)$,

\noindent
(II) $m_n((\xi_L)^{a},\eta,(\xi)^b,\th,(\xi_L)^c,\eta,(\xi)^d)\in\Ext^1(L,\OO)$,

\noindent
(III) $m_n((\xi)^{a},\th,(\xi_L)^b,\eta,(\xi)^c,\th,(\xi_L)^d,\eta,(\xi)^e)\in\Hom(\OO,\OO)$,

\noindent
(IV) $m_n((\xi_L)^{a},\eta,(\xi)^b,\th,(\xi_L)^c,\eta,(\xi)^d,\th,(\xi_L)^e)\in\Hom(L,L)$,

\noindent
where we denote by $(\xi)^a$ the string $(\xi,\ldots,\xi)$ with $\xi$ repeated $a$ times.

By the cyclic symmetry \eqref{cyclic-sym}, we have
$$\lan m_n((\xi_L)^{a},\eta,(\xi)^b,\th,(\xi_L)^c,\eta,(\xi)^d),\th\ran=
\lan m_n((\xi)^{b},\th,(\xi_L)^c,\eta,(\xi)^d,\th,(\xi_L)^a),\eta\ran,
$$
$$\lan m_n((\xi)^{a},\th,(\xi_L)^b,\eta,(\xi)^c,\th,(\xi_L)^d,\eta,(\xi)^e),\xi\ran=
\lan m_n((\xi)^{a+e+1},\th,(\xi_L)^b,\eta,(\xi)^c,\th,(\xi_L)^d),\eta\ran,$$
$$\lan m_n((\xi_L)^{a},\eta,(\xi)^b,\th,(\xi_L)^c,\eta,(\xi)^d,\th,(\xi_L)^e),\xi_L\ran=
\lan m_n((\xi)^{b},\th,(\xi_L)^c,\eta,(\xi)^d,\th,(\xi_L)^{a+e+1}),\eta\ran.$$
Hence, it is enough to compute the products of type (I), i.e., the
coefficients
$$\lan m_n((\xi)^{a},\th,(\xi_L)^b,\eta,(\xi)^c,\th,(\xi_L)^d),\eta\ran.$$

\subsection{Calculation I: combinatorial part}

We start by computing some signs $\eps(T)$.
For a pair of oriented planar rooted $3$-valent trees $T_1, T_2$ let us denote by $\join(T_1,T_2)$
the tree (of the same type) obtained by joining together the roots of $T_1$ and $T_2$ and
adding a root to the obtained new vertex (we keep $T_1$ on the left from $T_2$ in the plane).

\begin{lem}\label{sign-lem} 
Let $T=\join(T_1,T_2)$, where each $T_i$ has $n_i+1$ leaves ($i=1,2$), 
exactly one of which is marked by a degree $0$ element, and the rest marked by degree $1$ elements. 
Assume also that in each $T_i$ no two leaves of degree $1$ can be attached to the same vertex.
Then
$$\eps(T)=(-1)^{{n_1+n_2+2\choose 2}+n_2}.$$
\end{lem}

\Pf . By definition,
$$\eps(T)=\eps(T_1)\eps(T_2)\cdot (-1)^{(n_1+1)+n_1n_2},$$
so it remains to compute $\eps(T_i)$ for $i=1,2$.
Note that under our assumptions each tree $T_i$ has a very simple structure: it has the main stem
from the root to the leaf of degree $0$, to which leaves of degree $1$ can attach on the left and on the right:
\begin{figure}[htbp]
\begin{center}
\includegraphics[width=10cm]{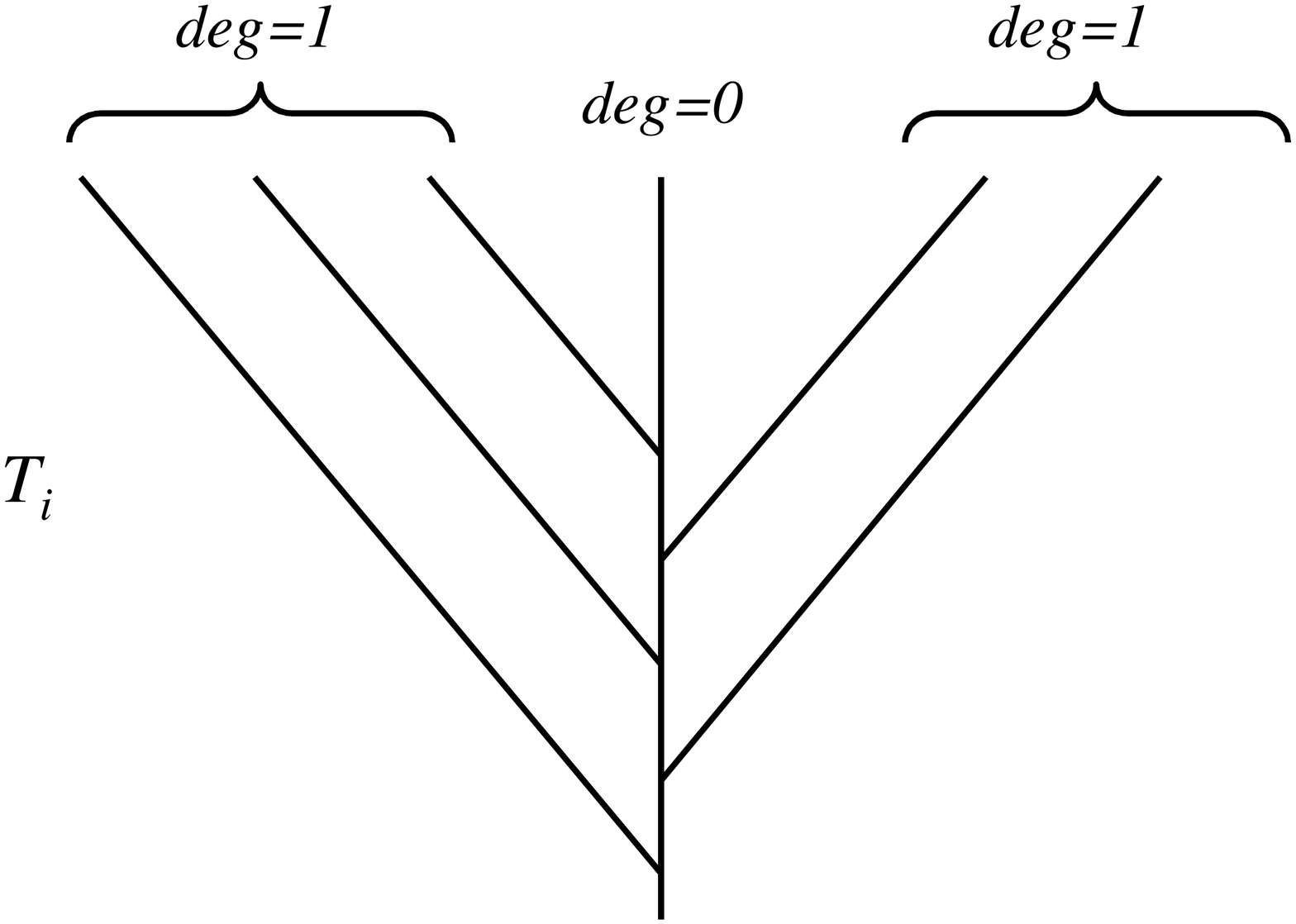}
\end{center}
\end{figure} \\
Suppose we have a pair of consecutive vertices $v$ and $w$ on the stem (with $v$ above $w$) such that
there is a degree $1$ leaf attaching to $v$ on the left and a degree $1$ leaf attaching to $w$ on the right.
Let $a$ be the number of degree $1$ leaves above $v$. Then the contribution of $v$ into the product
defining $\eps(T_i)$ equals $(-1)^{a+1}$, while the contribution of $w$ equals $(-1)^{a+2}$.
Hence, the contribution of both $v$ and $w$ is $-1$. It is easy to check that if $T'$ is the tree obtained
from $T$ by reversing the order of attaching these two leaves at vertices $v$ and $w$ then
the contribution of these vertices into $\eps(T')$ will still be $-1$. Assume that $T_i$
has exactly $a$ (resp., $b$) leaves to the left (resp., to the right) of the degree $0$ leaf,
and let $v_1,\ldots,v_a$ (resp., $w_1,\ldots,w_b$) be the vertices on the stem to which they attach.
By the above observation,
it is enough to consider the case when all the vertices $v_1,\ldots,v_a$ are above $w_1,\ldots,w_b$:
\begin{figure}[htbp]
\begin{center}
\includegraphics[width=10cm]{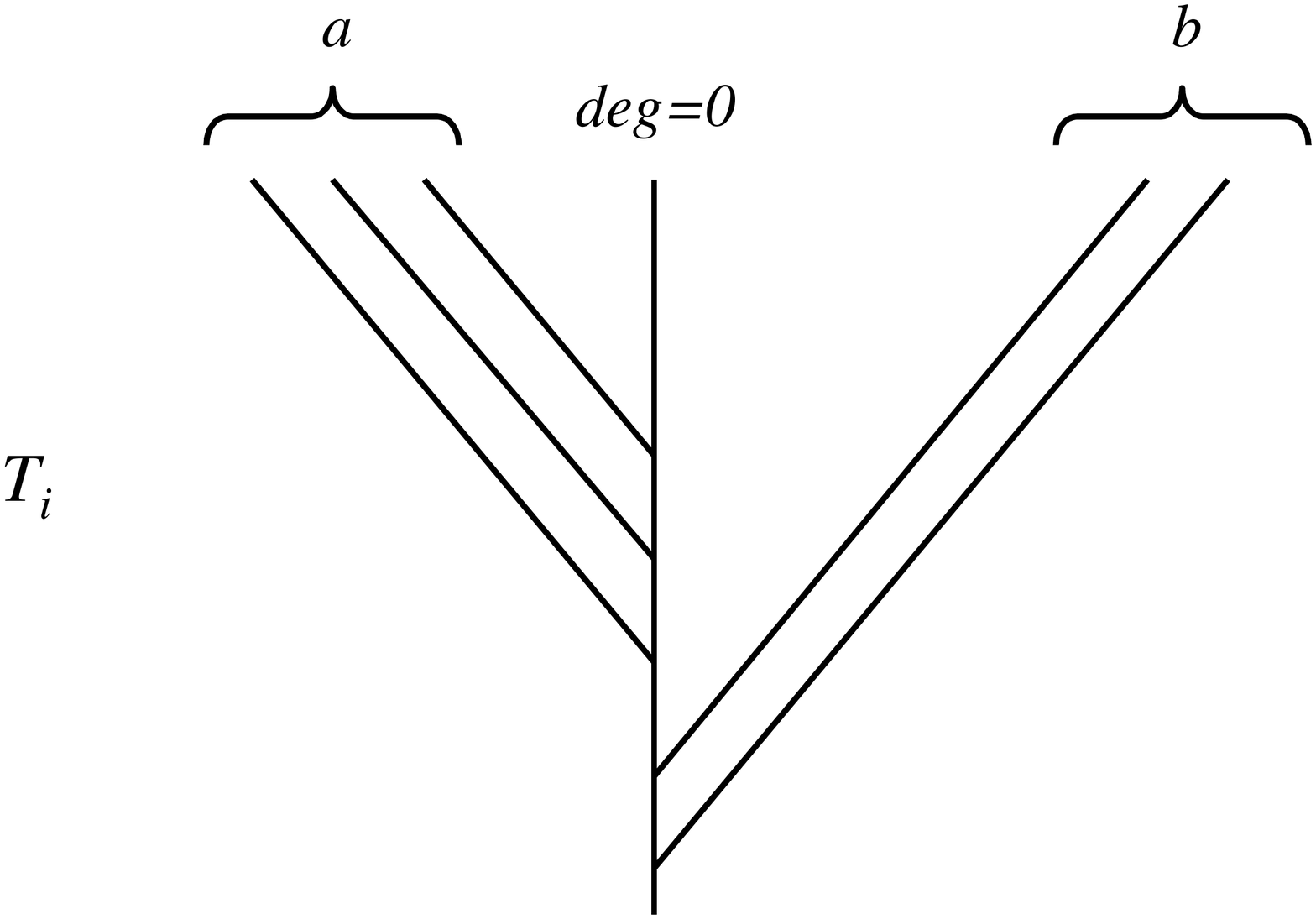}
\end{center}
\end{figure} \\
Then one can easily calculate that
$$\eps(T_i)=(-1)^{{a+b+1\choose 2}}.$$
Hence,
$$\eps(T)=(-1)^{{n_1+1\choose 2}+{n_2+1\choose 2}+(n_1+1)+n_1n_2}=(-1)^{{n_1+n_2+2\choose 2}+n_2}.$$
\ed

Next, we consider the terms $m_T((\xi)^{a},\th,(\xi_L)^b,\eta,(\xi)^c,\th,(\xi_L)^d)$.
Note that if two leaves of degree $1$ in $T$ are attached to the same vertex then $m_T$ vanishes because
it will involve taking products of two elements of degree $1$ in $A$. Henceforward, we assume
that no two leaves of degree $1$ can be attached to the same vertex.
It is convenient to introduce the operator 
$$H_L:C^{\infty}(L)\to C^{\infty}(L): s\mapsto Q(s\cdot d\ov{z})$$
and a similar operator $H_{\OO}$ on $C^{\infty}$-functions.
The next result follows immediately from the definition.

\begin{lem}\label{mT-lem} 
(i) Let $T=\join(T_1,T_2)$, where the leaves of $T_1$ are marked (from left to right) by
$$(\xi)^a,\th,(\xi_L)^b,\eta,(\xi)^{c_1},$$
while the leaves of $T_2$ are marked by
$$(\xi)^{c_2},\th,(\xi_L)^d.$$
Let $a=a_1+a_2$, where $a_1$ leftmost leaves of $T_1$ get attached to the main stem of $T_1$ below
the vertex where $\eta$ attaches to the main stem, and the next $a_2$ leaves get attached to the stem above this
vertex:
\begin{figure}[htbp]
\begin{center}
\includegraphics[width=10cm]{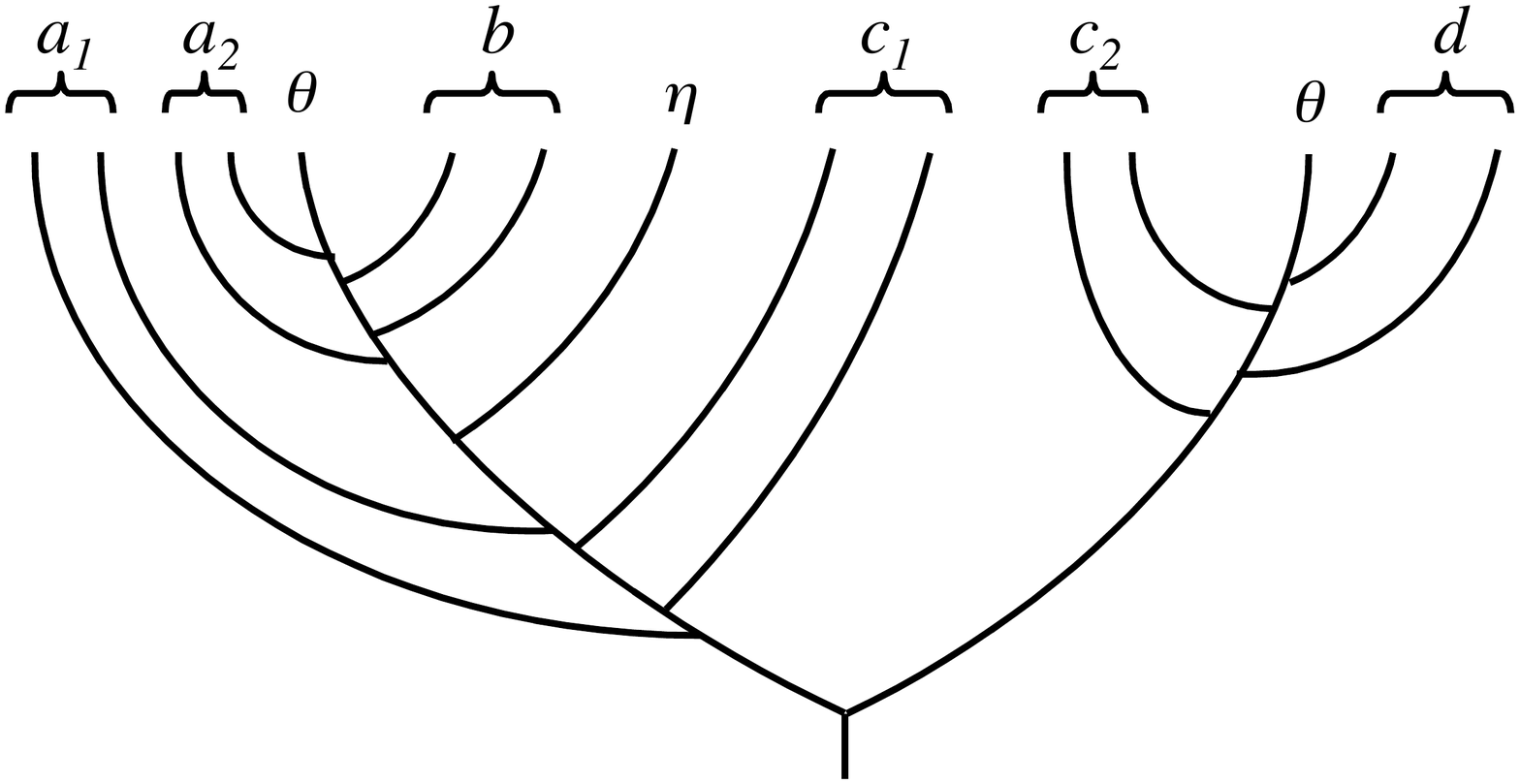}
\end{center}
\end{figure} \\
Then
$$m_T((\xi)^{a},\th,(\xi_L)^b,\eta,(\xi)^{c_1+c_2},\th,(\xi_L)^d)=\Pi
\left([H_{\OO}^{a_1+c_1}Q(H_L^{a_2+b}(\th)\cdot\eta)]\cdot H_L^{c_2+d}(\th)\right).$$

\noindent
(ii) Let $T=\join(T_1,T_2)$, where the leaves of $T_1$ are marked by
$$(\xi)^a,\th,(\xi_L)^{b_1},$$
while the leaves of $T_2$ are marked by
$$(\xi_L)^{b_2},\eta,(\xi)^c,\th,(\xi_L)^d.$$
Let $d=d_1+d_2$, where exactly $d_2$ rightmost leaves of $T_2$ get attached to the main stem of $T_2$ below
the vertex where $\eta$ attaches to the main stem. Then
$$m_T((\xi)^{a},\th,(\xi_L)^{b_1+b_2},\eta,(\xi)^c,\th,(\xi_L)^d)=\Pi
\left([H_{\OO}^{b_2+d_2}Q(H_L^{c+d_1}(\th)\cdot\eta)]\cdot H_L^{a+b_1}(\th)\right).$$
\end{lem}

Combining Lemmas \ref{sign-lem} and \ref{mT-lem} we arrive at the following expression for $m_n$.

\begin{lem}\label{comb-lem}
One has 
\begin{align*}
&\lan m_n((\xi)^{a},\th,(\xi_L)^b,\eta,(\xi)^c,\th,(\xi_L)^d),\eta\ran=\\
&(-1)^{{n\choose 2}+1}\cdot
\sum_{a=a_1+a_2;c=c_1+c_2}{a_2+b\choose a_2}{a_1+c_1\choose a_1}{c_2+d\choose c_2}
\phi(a_2+b,a_1+c_1,c_2+d)+\\
&(-1)^{{n\choose 2}+n+1}\cdot\sum_{b=b_1+b_2;d=d_1+d_2}{c+d_1\choose c}{b_2+d_2\choose b_2}{a+b_1\choose a}
\phi(c+d_1,b_2+d_2,a+b_1),
\end{align*}
where 
$$\phi(m,n,p)=(-1)^p\lan\Pi\left([H_{\OO}^nQ(H_L^m(\th)\cdot\eta)]\cdot H_L^p(\th)\right),\eta\ran.$$
\end{lem}

\subsection{Calculation II: analytic part}

We will use real coordinates $u,v$ on $\C$ such that $z=u+v\tau$.
Let us also set $a=\Im(\tau)$.
 
\begin{lem}\label{main-lem} Consider the differential operator $D=-\frac{a}{\pi}\frac{\pa}{\pa z}-2iav$.
Then for $k\ge 0$ one has:
\begin{equation}\label{HD-eq}
H_L^k\th=\frac{1}{k!}\cdot D^k\th=\frac{(-2ia)^k}{k!}\cdot\sum_{n\in\Z}(n+v)^k\exp(\pi i\tau n^2+2\pi i nz),
\end{equation}
\begin{equation}\label{Dtheta-eq}
\sqrt{2a}(D^k\th)\cdot \ov{\th}\exp(-2\pi a v^2)=\sum_{(m,n)\in\Z^2}(-1)^{mn}
(m\ov{\tau}-n)^k\exp(-\frac{\pi}{2a}|m\tau-n|^2+2\pi i (mu+nv)).
\end{equation}
\end{lem}

\Pf . The case $k=0$ of the identity \eqref{Dtheta-eq} is well-known 
(see \cite{P-hi-pr}, eq. (2.2)). The general case follows easily by applying $D^k$.
Next, let us prove \eqref{HD-eq} by induction in $k$.
Recall that the operator $Q:\Om^{0,1}(L)\to\Om^{0,0}(L)$ is uniquely determined by the following two properties:
$\dbar\circ Q=\id$ and the image of $Q$ is orthogonal to $\th$.
Thus, $H_L^k\th$ equals the unique function $f$ such that $\frac{\pa f}{\pa \ov{z}}=\frac{D^{k-1}\th}{(k-1)!}$ and
$(f,\th)=0$. Thus, it is enough to check the identities
\begin{equation}\label{Dk-exp}
D^k\th=(-2ia)^k\cdot\sum_{n\in\Z}(n+v)^k\exp(\pi i\tau n^2+2\pi i nz),
\end{equation}
\begin{equation}\label{dbarD}
\frac{\pa}{\pa\ov{z}}D^k\th=kD^{k-1}\th,
\end{equation}
$$(D^k\th,\th)=0$$
for $k>0$.
The latter follows immediately from the Fourier expansion
\eqref{Dtheta-eq}, since $(D^k\th,\th)$ is proportional to the Fourier coefficient of this expansion
corresponding to $(m,n)=(0,0)$. To check \eqref{Dk-exp} one can apply $D$ to the similar expansion
for $D^{k-1}\th$. Finally, one checks \eqref{dbarD} by applying $\frac{\pa}{\pa\ov{z}}$ to the right-hand side
of \eqref{Dk-exp}.
\ed

Now we can calculate the expressions $\phi(m,n,p)$ (see Lemma \ref{comb-lem}).

\begin{lem}\label{phi-lem}
One has 
$$\phi(k,l,p)=\frac{1}{k!p!}(\frac{a}{\pi})^{l+1}\cdot\sum_{\om\in\Z+\Z\tau\setminus\{0\}}
\frac{\ov{\om}^{k+p}}{\om^{l+1}}\exp(-\frac{\pi}{a}|\om|^2)=
\frac{1}{k!p!}(\frac{a}{\pi})^{k+l+p+1}f_{k+p,l+1}(\Z+\Z\tau).$$
\end{lem}

\Pf . It is easy to see that the operator $Q:\Om^{0,1}\to\Om^{0,0}$ is given by 
$$Q(\exp(2\pi i (mu+nv))d\ov{z})=\begin{cases}\frac{a}{\pi(m\tau-n)}\exp(2\pi i(mu+nv)), & (m,n)\neq (0,0),\\
0 & (m,n)=(0,0).
\end{cases}
$$
Hence, using Lemma \ref{main-lem} we obtain
\begin{align*}
&Q(H_L^k(\th)\cdot\eta)=\frac{1}{k!}Q(\sqrt{2a}(D^k\th)\cdot \ov{\th}\exp(-2\pi a v^2)=\\
&\frac{a}{k!\pi}\cdot
\sum_{(m,n)\in\Z^2\setminus\{(0,0)\}}(-1)^{mn}
\frac{(m\ov{\tau}-n)^k}{(m\tau-n)}\exp(-\frac{\pi}{2a}|m\tau-n|^2+2\pi i (mu+nv)).
\end{align*}
Therefore,
\begin{equation}\label{HOQHL-eq}
H_{\OO}^lQ(H_L^k(\th)\cdot\eta)=\frac{1}{k!}(\frac{a}{\pi})^{l+1}\cdot
\sum_{(m,n)\in\Z^2\setminus\{(0,0)\}}(-1)^{mn}
\frac{(m\ov{\tau}-n)^k}{(m\tau-n)^{l+1}}\exp(-\frac{\pi}{2a}|m\tau-n|^2+2\pi i (mu+nv)).
\end{equation}

Next, comparing the formulas for $\eta$ and for the metric on $L$ 
we observe that for a $C^{\infty}$-section $f$ of $L$ one has
$(f,\th)=0$ if and only if $\lan f,\eta\ran=0$. Hence, for $f\in\C^{\infty}(L)$ one has
$$\lan\Pi(f),\eta\ran=\lan f,\eta\ran$$
(since $\Pi$ is the orthogonal projection onto $\C\th$).
Therefore, 
$$\phi(k,l,p)=(-1)^p\lan [H^l_{\OO}Q(H_L^k(\th)\cdot\eta]\cdot H^p_L(\th),\eta\ran=
(-1)^p\lan H^l_{\OO}Q(H_L^k(\th)\cdot\eta, H^p_L(\th)\cdot\eta\ran.$$
Now the right-hand side can be computed using the Fourier expansion
for $H^p_L(\th)$ from Lemma \ref{main-lem} and the Fourier expansion \eqref{HOQHL-eq}:
$$\phi(k,l,p)=\frac{1}{k!p!}(\frac{a}{\pi})^{l+1}\cdot\sum_{(m,n)\in\Z^2\setminus\{(0,0)\}}
\frac{(m\ov{\tau}-n)^{k+p}}{(m\tau-n)^{l+1}}\exp(-\frac{\pi}{a}|m\tau-n|^2).$$
\ed

\subsection{Calculation III: conclusion}

It remains to put everything together. Substituting the expressions for $\phi(k,l,p)$ found in Lemma \ref{phi-lem}
into the formula of Lemma \ref{comb-lem} we get
\begin{align*}
&\lan m_n((\xi)^{a},\th,(\xi_L)^b,\eta,(\xi)^c,\th,(\xi_L)^d),\eta\ran=\\
&(-1)^{{n\choose 2}+1}(\frac{\Im(\tau)}{\pi})^{n-2}\cdot
\sum_{a=a_1+a_2;c=c_1+c_2}\frac{(a_1+c_1)!}{a_1!a_2!c_1!c_2!b!d!}f_{a_2+c_2+b+d,a_1+c_1+1}+\\
&(-1)^{{n\choose 2}+n+1}(\frac{\Im(\tau)}{\pi})^{n-2}\cdot
\sum_{b=b_1+b_2;d=d_1+d_2}\frac{(b_2+d_2)!}{b_1!b_2!d_1!d_2!a!c!}f_{b_1+d_1+a+c,b_2+d_2+1}.
\end{align*}
Since $f_{k,l}=0$ unless $k+l$ is even, this immediately implies that $m_n=0$ for odd $n$.
Now assuming that $n$ is even we can rewrite the above equation as follows (denoting $k=a_1+c_1$ and
$l=b_2+d_2$):
\begin{align*}
&(-1)^{{n\choose 2}+1}(\frac{\pi}{\Im(\tau)})^{n-2}\cdot \lan m_n((\xi)^{a},\th,(\xi_L)^b,\eta,(\xi)^c,\th,(\xi_L)^d),\eta\ran=\\
&\frac{1}{b!d!}\cdot\sum_{k\ge 0}C(a,c,k)f_{a+c-k+b+d,k+1}+
\frac{1}{a!c!}\cdot\sum_{l\ge 0}C(b,d,l)f_{b+d-l+a+c,l+1},
\end{align*}
where
$$C(a,c,k)=\sum_{a_1+c_1=k; a_1\le a, c_1\le c}\frac{k!}{a_1!c_1!(a-a_1)!(c-c_1)!}=
\frac{k!}{a!c!}\sum_{a_1+c_1=k}{a\choose a_1}{c\choose c_1}=\frac{k!}{a!c!}{a+c\choose k}$$
Thus, our formula for $\lan m_n((\xi)^{a},\th,(\xi_L)^b,\eta,(\xi)^c,\th,(\xi_L)^d),\eta\ran$ takes form
\begin{align*}
&(-1)^{{n\choose 2}+1}(\frac{\pi}{\Im(\tau)})^{n-2}\cdot \lan m_n((\xi)^{a},\th,(\xi_L)^b,\eta,(\xi)^c,\th,(\xi_L)^d),\eta\ran=\\
&\frac{1}{a!b!c!d!}\cdot\sum_{k\ge 0}k!\left({a+c\choose k}+{b+d\choose k}\right)f_{n-3-k,k+1}.
\end{align*}
Taking into account Corollary \ref{main-cor} we obtain for even $n$
\begin{equation}\label{main-eq}
\lan m_n((\xi)^{a},\th,(\xi_L)^{b},\eta, (\xi)^{c}, \th, (\xi_L)^{d}),\eta\ran=
(-1)^{{n\choose 2}+1}\frac{1}{a!b!c!d!}\cdot(\frac{\Im(\tau)}{\pi})^{n-2}\cdot g_{a+c,b+d}.
\end{equation}
Let us summarize our calculations. We set
\begin{equation}\label{Mabcd}
M(a,b,c,d):=(-1)^{{a+b+c+d+1\choose 2}}\frac{1}{a!b!c!d!}\cdot(\frac{\Im(\tau)}{\pi})^{a+b+c+d+1}\cdot g_{a+c,b+d}.
\end{equation}

\begin{thm}\label{main-thm} 
The only non-trivial higher products $m_n$ of the $A_{\infty}$-structure on $B=\Ext^*(\OO\oplus L,\OO\oplus L)$
are of the form
$$m_n((\xi)^{a},\th,(\xi_L)^{b},\eta,(\xi)^{c}, \th,(\xi_L)^{d})=M(a,b,c,d)\cdot\th,$$
$$m_n((\xi_L)^{a},\eta,(\xi)^{b},\th,(\xi_L)^{c}, \eta, (\xi)^{d})=M(a,b,c,d)\cdot\eta,$$
$$m_n((\xi)^{a},\th,(\xi_L)^{b},\eta,(\xi)^{c}, \th, (\xi_L)^{d},\eta,(\xi)^e)=M(a+e+1,b,c,d)\cdot\id_{\OO},$$
$$m_n((\xi_L)^{a},\eta,(\xi)^{b},\th,(\xi_L)^{c},\eta,(\xi)^{d},\th,(\xi_L)^e)=M(a+e+1,b,c,d)\cdot\id_L.$$
All products $m_n$ with odd $n$ vanish.
\end{thm}

\begin{rems} 1. Assume that $\tau$ belongs to a ring of integers of an imaginary quadratic field
(so that our elliptic curve admits complex multiplication). Set
$\varpi=2\pi|\eta(q)|^2$, where $\eta(q)$ is the Dedekind's $\eta$-function, and $q=\exp(2\pi i\tau)$.
Then the numbers
$\varpi^{-a-b-1}\cdot g_{a,b}$ are algebraic over $\Q$ (see \cite{Weil}, sec.~VI.6). 
Hence, if we multiply our basis elements of degree $1$ ($\eta$, $\xi$ and $\xi_L$) 
by the factor $|\eta(q)|^{-2}$ then the structure
constants of our $A_{\infty}$-structure with respect to the new basis will be algebraic over $\Q$.

\noindent 2. Another meaningful rescaling is obtained if we multiply our basis elements of degree
$1$ by the factor $\pi/\Im(\tau)$. Then the new structure constants $M'(a,b,c,d)$ will all have
limit at the cusp $\Im(\tau)\to +\infty$. Namely, using the well known $q$-series for the Eisenstein
series, one can easily check that the only nonzero limiting values will be 
$$M'(i,0,j,0)=M'(0,i,0,j)\to(-1)^{{i+j+1\choose 2}}\cdot {i+j\choose i}\cdot 2\zeta(i+j+1), 
$$
where $i+j$ is odd.
\end{rems}

\subsection{$A_{\infty}$-constraint}

The $A_{\infty}$-axiom (we follow \cite{Keller}, 3.1 for sign conventions) gives certain quadratic equations on coefficients 
$(M(a,b,c,d))$ and hence leads to identities for $g_{m,n}$.
For example, applying this axiom to the string 
$$(\xi)^a,\th,(\xi_L)^b,\eta,(\xi)^c,\th,(\xi_L)^d,\eta,(\xi)^e,\th,(\xi_L)^f,$$
where $a,b,c,d,e,f$ are positive, we get 
\begin{align*}
&\sum_{a=a_1+a_2;d=d_1+d_2}(-1)^{(a_2+b+c+d_1+1)(a_1+d_2+e+f)+a_1}M(a_2,b,c,d_1)M(a_1,d_2,e,f)+\\
&\sum_{b=b_1+b_2;e=e_1+e_2}(-1)^{(b_2+c+d+e_1+1)(a+b_1+e_2+f+1)+a+b_1+1}M(b_2,c,d,e_1)M(a,b_1,e_2,f)+\\
&\sum_{c=c_1+c_2;f=f_1+f_2}(-1)^{(c_2+d+e+f_1+1)(a+b+c_1+1+f_2)+a+b+c_1}M(c_2,d,e,f_1)M(a,b,c_1,f_2)=0.
\end{align*}
When one of $a,b,c,d,e,f$ is zero, additional terms will arise due to the presence of double products.
For example, for the string
$$(\xi)^a,\th,\eta,\th,\eta,\th,(\xi_L)^b$$
we get
\begin{align*}
&\sum_{a=a_1+a_2}(-1)^{(a_2+1)(a_1+b)+a_1}M(a_2,0,0,0)M(a_1,0,0,b)+\\
&\sum_{b=b_1+b_2}(-1)^{(b_1+1)(a+b_2+1)+a}M(0,0,0,b_1)M(a,0,0,b_2)+\\
&(-1)^a[M(a+1,0,0,b)-M(a,1,0,b)+M(a,0,1,b)-M(a,0,0,b+1)]+\\
&\de_{b,0}(-1)^aM(a+1,0,0,0)-\de_{a,0}M(b+1,0,0,0)=0.
\end{align*}
Substituting into the above identities the expressions for $M(a,b,c,d)$ from \eqref{Mabcd} we arrive
at the following result.

\begin{prop}\label{id-prop} 
(i) For positive integers $a,b,c,d,e,f$ one has
\begin{align*}
&\sum_{a=a_1+a_2;d=d_1+d_2}(-1)^{c+d_1+1}{a\choose a_1}{d\choose d_1}g_{a_2+c,b+d_1}
g_{a_1+e,d_2+f}+\\
&\sum_{b=b_1+b_2;e=e_1+e_2}(-1)^{b_2+1}{b\choose b_1}{e\choose e_1}g_{b_2+d,c+e_1}
g_{a+e_2,b_1+f}+\\
&\sum_{c=c_1+c_2;f=f_1+f_2}(-1)^{c_1}{c\choose c_1}{f\choose f_1}g_{c_2+e,d+f_1}
g_{a+c_1,b+f_2}=0.
\end{align*}

\noindent
(ii) For integers $a,b\ge 0$ one has
$$\sum_{a=a_1+a_2}{a\choose a_1}g_{a_1,0}g_{a_2,b}-\frac{a+2+\de_{b,0}}{a+1}g_{a+1,b}=
\sum_{b=b_1+b_2}{b\choose b_1}g_{0,b_1}g_{a,b_2}-\frac{b+2+\de_{a,0}}{b+1}g_{a,b+1}.
$$
\end{prop}

Note that the identity in (ii) gives a recursive formula for $g_{a+1,b}$ in terms of all $g_{a',b'}$ with $a'\le a$.
Since $g_{0,n}=n!e^*_{n+1}$, we recover the fact that all $g_{a,b}$ are polynomials in $(e^*_n)$ with rational
coefficients. For example, in the case $a=0$ the obtained identity (for even $n$)
$$2g_{1,n}=-\sum_{n=m+k}{n\choose m}g_{0,m}g_{0,k}+\frac{n+3}{n+1}g_{n+1,0}$$
is equivalent to the formula
$$\frac{1}{n}We^*_n=-\sum_{n=m+k}e^*_{m+1}e^*_{k+1}+(n+3)e_{n+2}$$
found in \cite{Weil}, VI.5.

\end{document}